\documentclass[12pt,a4paper]{article}
\usepackage{amsmath}
\usepackage{amsthm}
\usepackage{amsfonts}
\usepackage{amssymb}
\usepackage{latexsym}
\usepackage[usenames]{color}

\newtheorem{theorem}{Theorem}
\newtheorem{lemma}{Lemma}
\newtheorem{corollary}{Corollary}
\newtheorem{problem}{Problem}

\theoremstyle{definition}
\newtheorem{definition}{Definition}

\newtheorem{example}{Example} 

\begin{document}

\title{\Large{\textbf{Equational quasigroup  definitions}}}
\author{\normalsize {V.A.~Shcherbacov$^{\ast}$, D.I. Pushkashu$^{\ast}$, A.V.~Shcherbacov$^{\ast\ast}$}
}

 \maketitle

\begin{abstract}
\noindent Quasigroup equational definitions are given.

\medskip

\noindent \textbf{2000 Mathematics Subject Classification:} 20N05

\medskip

\noindent \textbf{Key words and phrases:} quasigroup, equational quasigroup
\end{abstract}


\bigskip


Basic and standard definition of a binary quasigroup is the following
\begin{definition} \label{MAIN_QUAS_DEF} Binary groupoid $(Q, \circ)$ is called a quasigroup if for all ordered pairs $(a, b)\in Q^2$
there exist unique solutions $x, y \in Q$ to the equations $x \circ a = b$ and $a \circ y = b$ \cite{VD}.
\end{definition}
\index{quasigroup!binary}

 T.~Evans has given equational definition of a quasigroup \cite{EVANS} (see  Definition
\ref{def3}). Evans' definition usually used by the study of universal algebraic questions of quasigroup theory.
The equivalence of Definitions \ref{MAIN_QUAS_DEF} and \ref{def3} is well known fact \cite{VD, SCERB_03}.

In this paper we give some new   equational  definitions of quasigroups.
The subject of this paper  is close with the subject of  articles \cite{SCERB_07, SCERB_TAB_PUSH_09}.

We shall use basic terms and concepts from books \cite{VD, 1a, HOP}.

Garrett Birkhoff in his  famous  book \cite{BIRKHOFF}   defined equational quasigroup as a algebra with three binary operations $(Q, \cdot, \slash, \backslash)$ that fulfil the following six identities
\begin{equation}
x\cdot(x \backslash y) = y \label{A}
\end{equation}
\begin{equation}
(y / x)\cdot x = y \label{C}
\end{equation}
\begin{equation}
x\backslash (x \cdot y) = y \label{B}
\end{equation}
\begin{equation}
(y \cdot x)/ x = y \label{D}
\end{equation}
\begin{equation}
x/(y\backslash x) = y  \label{T}
\end{equation}
\begin{equation}
(x/y)\backslash x = y \label{R}
\end{equation}

Results of this paper are connected with the following
\begin{problem}
Research properties of algebra $(Q, \cdot, \slash, \backslash)$ with various combinations
of identities (\ref{A})--(\ref{R})  (\cite{SCERB_03}, page 11).
\end{problem}

It is well known the following

\begin{lemma}
In algebra  $(Q, \cdot, \backslash, /)$ with identities (\ref{A}), (\ref{C}), (\ref{B}), and (\ref{D}) identities (\ref{T}) and (\ref{R}) are true   \cite{JDH_2007, SCERB_03, SCERB_07}.
\end{lemma}
\begin{proof}
We can re-write  identity (\ref{D}) in the following form
\begin{equation} \label{four}
(x \cdot (x\backslash y)) \slash (x\backslash y) = x
\end{equation}
By identity (\ref{A})  $x \cdot (x\backslash y) = y$. Thus  from identity  (\ref{four}) we obtain
$y \slash (x\backslash y) = x$, i.e. we obtain identity (\ref{T}).

We can re-write  identity (\ref{B}) in the following form
\begin{equation} \label{twor}
(x\slash y) \backslash ((x\slash y)\cdot y) = y
\end{equation}
By identity (\ref{C})  $(x\slash y)\cdot y = x$. Thus  from identity  (\ref{twor}) we obtain
$(x\slash y) \backslash x = y$, i.e. we obtain identity (\ref{R}).
\end{proof}

Therefore  it is used the following T.~Evans'  equational definition of a quasigroup   \cite{EVANS}.
\begin{definition}    An algebra $(Q, \cdot, \backslash, /)$ with identities (\ref{A}),
(\ref{C}), (\ref{B}) and (\ref{D}) is called a \textit{quasigroup} \cite{EVANS, BIRKHOFF, VD, 1a,
HOP, BURRIS}. \label{def3}
\end{definition}

\begin{lemma}  \label{CT_A}  In  algebra $(Q, \cdot, \backslash, /)$ from identities (\ref{C}) and (\ref{T}) it follows identity (\ref{A}).
\end{lemma}
\begin{proof}
We can re-write  identity (\ref{C}) in the following form
\begin{equation} \label{seven}
(x \slash (y\backslash x)) \cdot (y\backslash x) = x
\end{equation}
But by identity (\ref{T}) $x \slash (y\backslash x)=y$. Therefore we can rewrite identity (\ref{seven})
in the following form
\begin{equation} \label{eiht}
y \cdot (y\backslash x) = x
\end{equation}
Then  we obtain identity (\ref{A}).
\end{proof}

\begin{lemma}  \label{BT_D}  In  algebra $(Q, \cdot, \backslash, /)$ from identities (\ref{B}) and (\ref{T}) it follows identity (\ref{D}).
\end{lemma}
\begin{proof}
We can re-write  identity (\ref{T}) in the following form
\begin{equation} \label{nine}
(x \cdot y) \slash (x \backslash (x\cdot y)) = x
\end{equation}
But by identity (\ref{B}) $x \backslash (x\cdot y)= y$. Therefore identity (\ref{nine}) takes the form
$(x \cdot y) \slash y = x$ and it coincides with identity (\ref{D}).
\end{proof}

\begin{lemma}  \label{RD_B}  In  algebra $(Q, \cdot, \backslash, /)$ from identities (\ref{D}) and (\ref{R}) it follows identity (\ref{B}).
\end{lemma}
\begin{proof}
We can re-write  identity (\ref{R}) in the following form
\begin{equation} \label{ten}
((x \cdot y) \slash y) \backslash (x\cdot y) = y
\end{equation}
But by identity (\ref{D}) $(x \cdot y) \slash y= x$. Therefore identity (\ref{ten}) takes the form
$x  \backslash (x \cdot y) = y$ and it coincides with identity (\ref{B}).
\end{proof}

\begin{lemma}  \label{AR_C}  In  algebra $(Q, \cdot, \backslash, /)$ from identities (\ref{A}) and (\ref{R}) it follows identity (\ref{C}).
\end{lemma}
\begin{proof}
We can re-write  identity (\ref{A}) in the following form
\begin{equation} \label{eleven}
(x \slash y) \cdot ((x \slash y) \backslash x) = x
\end{equation}
But by identity (\ref{R}) $(x \slash y) \backslash x = y$. Therefore identity (\ref{eleven}) takes the form
$(x \slash y) \cdot y  = x$ and it coincides with identity (\ref{C}).
\end{proof}

\begin{theorem}  \label{CBT_TH}  An algebra $(Q, \cdot, \backslash, /)$ with identities
(\ref{C}), (\ref{B}) and (\ref{T}) is  a \textit{quasigroup}.
\end{theorem}
\begin{proof}
The proof follows from Lemmas \ref{CT_A} and \ref{BT_D}. \end{proof}

\begin{theorem}  \label{ADR_TH}  An algebra $(Q, \cdot, \backslash, /)$ with identities
(\ref{A}), (\ref{D}) and (\ref{R}) is  a \textit{quasigroup}.
\end{theorem}
\begin{proof}
The proof follows from Lemmas \ref{RD_B} and \ref{AR_C}.
\end{proof}

In the following corollary  we give definitions of equational quasigroup using four identities from the identities (\ref{A})--(\ref{R}).

\begin{corollary}  \label{FOUR_ID_ADR_TH}
\begin{enumerate}
  \item  An algebra $(Q, \cdot, \backslash, /)$ with identities
(\ref{A}), (\ref{C}), (\ref{B}) and (\ref{T}) is  a \textit{quasigroup}.
  \item An algebra $(Q, \cdot, \backslash, /)$ with identities
(\ref{C}), (\ref{B}), (\ref{D}) and (\ref{T}) is  a \textit{quasigroup}.
  \item An algebra $(Q, \cdot, \backslash, /)$ with identities
(\ref{A}), (\ref{C}), (\ref{D}) and (\ref{R}) is  a \textit{quasigroup}.
  \item An algebra $(Q, \cdot, \backslash, /)$ with identities
(\ref{A}), (\ref{B}), (\ref{D}) and (\ref{R}) is  a \textit{quasigroup}.
  \item An algebra $(Q, \cdot, \backslash, /)$ with identities
(\ref{A}), (\ref{B}), (\ref{T}) and (\ref{R}) is  a \textit{quasigroup}.
  \item An algebra $(Q, \cdot, \backslash, /)$ with identities
(\ref{A}), (\ref{D}), (\ref{T}) and (\ref{R}) is  a \textit{quasigroup}.
  \item An algebra $(Q, \cdot, \backslash, /)$ with identities
(\ref{C}), (\ref{B}), (\ref{T}) and (\ref{R}) is  a \textit{quasigroup}.
\item An algebra $(Q, \cdot, \backslash, /)$ with identities
(\ref{C}), (\ref{D}), (\ref{T}) and (\ref{R}) is  a \textit{quasigroup}.
\end{enumerate}
\end{corollary}
\begin{proof} The proof follows from Theorems \ref{CBT_TH} and  \ref{ADR_TH},  Lemmas \ref{CT_A}, \ref{BT_D},   \ref{RD_B} and \ref{AR_C}.
\end{proof}

The  proofs of Lemmas \ref{CT_A}, \ref{BT_D},   \ref{RD_B} and \ref{AR_C} are obtained using Prover 9 \cite{MAC_CUNE_PROV}.

Information on properties of algebras with 3-element sets  of identities that are taken from identities (\ref{A})--(\ref{D}) it is possible to deduce   from results of the articles \cite{SCERB_07, SCERB_TAB_PUSH_09}.
For example, algebra $(Q, \cdot, \slash, \backslash)$ with identities (\ref{A}), (\ref{C}), (\ref{B}) is  a left quasigroup with right division.

\begin{example} \label{EXAMPLE left_quas_DIVIZ_RIGHT} Let  $x\circ y =  \left[ x/2 \right] - 1\cdot y $ for all
$x, y \in {\mathbb Z}$, where $ ({\mathbb Z}, +, \cdot) $ is the ring of integers, $\left[ x/2 \right] = a,$ if $x=2\cdot a$;
$\left[ x/2 \right] = a$, if $x=2\cdot a + 1$. It  is possible to check that $ ({\mathbb Z}, \circ) $ is a left quasigroup with right division.
\end{example}

\noindent \footnotesize
{$^{\ast}$Institute of Mathematics and Computer Science \hfill $^{\ast\ast}$Moldova State University\\
Academy of Sciences of Moldova \hfill A. Mateevici str. 60, MD-2009, Chi\c sin\u au\\
5 Academiei str., MD$-$2028, Chi\c{s}in\u{a}u \hfill  Moldova\\
 Moldova  \hfill E-mail: \emph{admin@sibirsky.org }\\
E-mail: \emph{scerb@math.md } \\ \emph{dmitry.pushkashu@gmail.com}}

\end{document}